\documentclass[a4paper,reqno]{amsart}
\usepackage{geometry}
\usepackage{mathrsfs}
\usepackage{syntonly}
\usepackage{amsmath}
\usepackage{amsthm}
\usepackage{amsfonts}
\usepackage{amssymb}
\usepackage{latexsym}
\usepackage{amscd,amssymb,amsopn,amsmath,amsthm,graphics,amsfonts,mathrsfs,accents,enumerate,verbatim,calc}
\usepackage[dvips]{graphicx}
\usepackage[colorlinks=true,linkcolor=blue,citecolor=blue]{hyperref}
\input xy
\xyoption{all}

\usepackage{color}

\numberwithin{equation}{section}

\theoremstyle{plain}

\newtheorem{thm}{Theorem}[section]
\newtheorem{cor}[thm]{Corollary}
\newtheorem{lem}[thm]{Lemma}
\newtheorem{prop}[thm]{Proposition}

\newtheorem{defn}[thm]{Definition}

\newtheorem{rem}[thm]{Remark}

\newdir{ >}{{}*!/-5pt/\dir{>}}

\newcommand{\Mod}{\operatorname{Mod}\nolimits}

\newcommand{\add}{\operatorname{add}\nolimits}
\newcommand{\Fac}{\operatorname{Fac}\nolimits}

\newcommand{\Hom}{\operatorname{Hom}\nolimits}

\newcommand{\Ext}{\operatorname{Ext}\nolimits}

\newcommand{\op}{\operatorname{op}\nolimits}

\renewcommand{\Mod}{\mathsf{Mod}\hspace{.01in}}
\renewcommand{\mod}{\mathsf{mod}\hspace{.01in}}

\newcommand{\M}{\mathcal M}

\newcommand{\B}{\mathcal B}

\newcommand{\U}{\mathcal U}
\newcommand{\V}{\mathcal V}
\newcommand{\A}{\mathcal A}

\newcommand{\h}{\mathcal H}
\newcommand{\s}{\mathcal S}
\newcommand{\T}{\mathcal T}

\newcommand{\N}{\mathcal N}
\newcommand{\R}{\mathcal R}
\newcommand{\X}{\mathscr X}
\newcommand{\Y}{\mathcal Y}
\newcommand{\Z}{\mathcal Z}
\newcommand{\C}{\mathcal C}

\renewcommand{\emph}{\textit}
\renewcommand{\phi}{\varphi}

\begin{document}

\title{$\tau$-tilting theory in abelian categories}\footnote{Yu Liu was supported by the Fundamental Research Funds for the Central Universities (Grant No. 2682018ZT25) and the National Natural Science Foundation of China (Grant No. 11901479). Panyue Zhou was supported by the National Natural Science Foundation of China (Grant Nos. 11901190 and 11671221), and by the Hunan Provincial Natural Science Foundation of China (Grant No. 2018JJ3205), and by the Scientific Research Fund of Hunan Provincial Education Department (Grant No. 19B239).}
\author{Yu Liu and Panyue Zhou}
\address{School of Mathematics, Southwest Jiaotong University, 610031, Chengdu, Sichuan, People's Republic of China}
\email{liuyu86@swjtu.edu.cn}
\address{College of Mathematics, Hunan Institute of Science and Technology, 414006, Yueyang, Hunan, People's Republic of China}
\email{panyuezhou@163.com}
\thanks{The authors would like to thank Professor Dong Yang and Professor Bin Zhu for helpful discussions.}
\begin{abstract}
Let $\A$ be a Hom-finite abelian category with enough projectives. In this note,
we show that any covariantly finite $\tau$-rigid subcategory is contained in a support $\tau$-tilting subcategory. We also show that
support $\tau$-tilting subcategories are in bijection with certain finitely generated torsion classes. Some applications of our main results are also given.
\end{abstract}
\maketitle

\section{Introduction}

In mathematics, especially representation theory, classical tilting theory describes a way to relate the module categories of two algebras using so-called tilting modules and associated tilting functors. Classical tilting theory was motivated by the reflection functors introduced by Bernstein, Gelfand and  Ponomarev \cite{BGP}. These functors were reformulated by Auslander, Platzeck and Reiten \cite{APR}, and generalized by Brenner and Butler \cite{BB}. It is well known that any almost complete tilting module
has either one or two complements, and mutation is possible only when there are two complements.
To make mutation always possible, it is desirable to enlarge
our class of tilting modules to get a more regular property that
almost complete ones always have two complements.  This is
accomplished by introducing $\tau$-tilting modules, or more precisely,
support $\tau$-tilting modules. It is called $\tau$-tilting theory by Adachi, Iyama and Reiten \cite{AIR}, which  generalizes classical tilting theory.
$\tau$-tilting modules can be viewed as a generalization of tilting modules.
They proved that any $\tau$-rigid module is a direct summand of some $\tau$-tilting module.
In fact, it is shown in \cite[Theorem 2.18]{AIR} that any almost complete support $\tau$-tilting module is the direct summand of exactly two support $\tau$-tilting modules. This means that mutation of support $\tau$-tilting modules is always possible.
\smallskip

Recently, Iyama, J{\o}rgensen and Yang \cite{IJY} gave a functor version of $\tau$-tilting theory.
They showed that two-term silting subcategories are in bijection with
support $\tau$-tilting subcategories under certain assumptions. They also showed that support $\tau$-tilting subcategories are in bijection with certain finitely generated torsion classes.
Let $\T$ be a triangulated category with shift functor $[1]$ and $\mathcal R$ be a rigid subcategory of
$\T$. Zhou and Zhu \cite{ZZ} proved that there exists a bijection between the set of two-term weak $\mathcal R[1]$-cluster tilting subcategories of $\T$ and the set of support $\tau$-tilting subcategories of $\mod\R$.
When $\mathcal R$ is a silting subcategory, the bijection above induces the bijection given by Iyama, J{\o}rgensen and Yang \cite{IJY}.
\smallskip

In this article, we discuss an abelian version of $\tau$-tilting theory. Our first main result is the following completion for $\tau$-rigid subcategories.

\begin{thm}{\rm (see Theorem \ref{main1} for details)}
Let $\A$ be an abelian category with enough projectives. Then any covariantly finite $\tau$-rigid subcategory is contained in a support $\tau$-tilting subcategory.
\end{thm}

Iyama, J{\o}rgensen and Yang \cite{IJY} also investigated the relation between $\tau$-tilting theory and torsion classes. More precisely,
they showed that support $\tau$-tilting pairs correspond bijectively with certain finitely generated torsion classes.
For a subcategory $\M$ of an abelian category, denote by $\Fac\M$ the subcategory consisting of factor
objects of finite direct sums of objects of $\M$, and $\mathbf{P}(\M)$ the Ext-projective objects of $\M$. Based on their ideas, our second main result is the following.

\begin{thm}{\rm (see Theorem \ref{main2} for details)}
Let $\A$ be an abelian category with enough projectives. There exists a bijection $\M\mapsto\Fac\M$ from the first of the following sets to the second:
\begin{itemize}
\item[(a)] Support $\tau$-tilting subcategories $\M$.
\smallskip

\item[(b)] Subcategories $\s=\Fac \mathbf{P}(\s)$ such that every projective object admits a left $\mathbf{P}(\s)$-approximation.
\end{itemize}
\end{thm}

This article is organized as follows. In Section 2, we give some basic concepts which be used later.
In Section 3, we prove our two main results. In Section 4, We give some applications of our main results.
\medskip

\section{Preliminaries}

In this article, let $k$ be a field and $\A$ be a Hom-finite abelian category over $k$. We assume that $\A$ has enough projectives $\mathcal P$.  For a subcategory $\B$, denote $\Fac \B$ by the subcategory
$$\{C\in \A \text{ }|\text{ there exists an epimorphism } B\to C\to 0, \mbox{where}~ B\in \B \};$$
denote  $\rm\bf P(\B)$ by the subcategory
$$\{D\in  \B \text{ }|\text{ } \Ext^1_{\A}(D,B)=0 ~\mbox{ for any } B\in \B \}.$$
Such object $D$ is called \emph{Ext-projective}.
Since we assume that $\A$ has enough projectives $\mathcal P$, we have $\A\simeq \mod \mathcal P$, see \cite[Corollaries 3.9 and 3.10]{B} for more details.
\smallskip

%The following definition is a special case of \cite[Definition 1.3]{IJY}.
Based on \cite[Definition 1.3]{IJY}, we give the following definition.

\begin{defn}\label{tau}
Let $\A$ be an abelian category with enough projectives $\mathcal P$.
\begin{itemize}
\item[(i)] A subcategory $\M$ of $\A$ is said to be $\tau$-rigid if any object $M\in \M$ admits an exact sequence $P_1\xrightarrow{f} P_0\to M\to 0$ such that $P_1,P_0\in \mathcal P$ and $\Hom_{\A}(f, M')$ is a surjection for any $M'\in \M$.
\smallskip

\item[(ii)] A $\tau$-rigid subcategory $\M$ is said to be support $\tau$-tilting if any projective object $P$ admits an exact sequence $P\xrightarrow{m} M^0\to M^1\to 0$ such that $M^0,M^1\in \M$ and $m$ is a left $\M$-approximation.

\item[(iii)] A $\tau$-rigid subcategory $\M$ is said to be $\tau$-tilting if any projective object $P$ admits an exact sequence $P\xrightarrow{m} M^0\to M^1\to 0$ such that $M^0,M^1\in \M$ and $m$ is a non-zero left $\M$-approximation.
\end{itemize}
\end{defn}

%\begin{rem}
%By the definition of a support $\tau$-tilting subcategory, the morphism $m$ in Definition \ref{tau} is a left $(\Fac \M)$-approximation.
%\end{rem}

Since $\mod \mathcal P$ is a full abelian subcategory of $\Mod \mathcal P$, and any monomorphism (resp. epimorphism, kernel, cokernel) in $\mod \mathcal P$ is also a monomorphism (resp. epimorphism, kernel, cokernel) in $\Mod \mathcal P$, by the proof of \cite[Lemma 5.2]{IJY}, we get the following important lemma.

\begin{lem}\label{imp}
A subcategory $\M\subseteq \A$ is $\tau$-rigid if and only if $\Ext^1_{\A}(\M,\Fac \M)=0$.
\end{lem}

\begin{rem}
$\Fac \M$ is called a finitely generated torsion class if it is closed under extensions. By the proof of \cite[Proposition 5.3]{IJY}, if $\M$ is a $\tau$-rigid subcategory, then $\Fac \M$ is a finitely generated torsion class.
\end{rem}

Let $\T$ be a Krull-Schimdt, Hom-finite triangulated category over $k$ with the shift functor $[1]$ and $\R$ be a rigid subcategory. Let
$$\h=\{T\in\T \text{ }|\text{ }\text{there exists a triangle }R^1\to T \to R^0[1]\to R^1[1], ~\mbox{where}~ R^0,R^1\in \R \}.$$
By the definition, We can get the following useful lemmas.

\begin{lem}\label{summand}  {\rm \cite[Proposition 2.1]{IY}}
$\h$ is closed under direct summands.
\end{lem}

\begin{lem}\label{inH}
If $T\in \T$ admits a triangle $R\to S\to T\to R[1]$ where $S\in\h$ and $R\in \R$, then $T\in \h$.
\end{lem}

For objects $A,B\in\T$ and a subcategory $\T'$ of $\T$, let $[\T'](A,B)$ be the subgroup of $\Hom_{\T}(A,B)$ consisting of morphisms which factor through objects in $\T'$. For a subcategory $\T_0$, we denote by $\overline{\T_0}$ the category which has the same objects as $\T_0$, and
$$\Hom_{\overline{\T_0}}(A,B)=\Hom_{\T_0}(A,B)/[\R[1]](A,B)$$
where $A,B\in \T_0$. For any morphism $f\in \Hom_{\T_0}(A,B)$, we denote its image in $\Hom_{\overline{\T_0}}(A,B)$ by $\overline f$. We know that $\overline{\T_0}$ is a subcategory of $\overline{\T}$.

\begin{defn}
Let $\U$ be a subcategory of $\T$.
\begin{itemize}
\item[(a)] $\U$ is called two-term $\R[1]$-rigid if $\U\subseteq \h$ and $[\R[1]](\U,\U[1])=0$.
\item[(b)] $\U$ is called two-term weak $\R[1]$-cluster tilting if $\U\subseteq \h$ and
$$\U=\{T\in \h \text{ }|\text{ } [\R[1]](T,\U[1])=0~ \text{ ~and~ } ~[\R[1]](\U,T[1])=0  \}.$$
\end{itemize}
\end{defn}

By \cite[Corollary B]{LC}, the bounded derived category $\mathrm{D}^b(\A)$ is a Krull-Schmidt, Hom-finite triangulated category over $k$. Since $\A$ is a full abelian subcategory (which is the heart of a $t$-structure) in $\mathrm{D}^b(\A)$, it is also Krull-Schmidt. We also assume $\A$ has enough projectives $\mathcal P$, then $\mathcal P$ is a rigid subcategory of $\mathrm{D}^b(\A)$.
We have a functor:
\begin{eqnarray*}
\mathbb{F}\colon \mathrm{D}^b(\A) & \longrightarrow & \Mod \mathcal P\\
 M& \longmapsto & \Hom_{\mathrm{D}^b(\A)}(-,M)|_{\mathcal P}
\end{eqnarray*}
Let $\T=\mathrm{D}^b(\A)$ and $\R=\mathcal P$. Then we can use the previous notions. For a morphism $f$ in $\h$, $\mathbb{F}|_{\h}(f)=0$ if and only if $f$ factors through $\mathcal P[1]$. $\mathbb{F}|_{\h}$ also induces an equivalence $H:\overline \h\to \mod \mathcal P$. Denote its quasi-inverse by $H^{-1}$. Since $\Hom_{\mathrm{D}^b(\A)}(-,-)|_{\A\times \A^{\op}}=\Hom_{\A}(-,-)$, we have an equivalence $\mathbb{F}|_{\A}\colon \A\to \mod \mathcal P$.
\smallskip

In the next section, we will use the following result frequently, which is a special case of \cite[Theorem 4.4, Theorem 4.5]{ZZ}.

\begin{thm}\label{ZZ}
The functor $\mathbb{F}\colon \mathrm{D}^b(\A)  \longrightarrow  \Mod \mathcal P$ induces a bijection $\Phi: \X\to \mathbb{F}(\X)$ from the first of the following sets to the second:
\begin{itemize}
\item[(I)] Two-term $\mathcal P[1]$-rigid subcategories of $\mathrm{D}^b(\A)$.
\smallskip

\item[(II)] $\tau$-rigid subcategories of $\mod \mathcal P$.
\end{itemize}
It also induces a bijection from the first of the following sets to the second:
\begin{itemize}
\item[(I)] Two-term weak $\mathcal P[1]$-cluster tilting subcategories of $\mathrm{D}^b(\A)$.
\smallskip

\item[(II)] Support $\tau$-tilting subcategories of $\mod \mathcal P$.
\end{itemize}
\end{thm}

\section{$\tau$-rigid subcategories and support $\tau$-tilting subcategories}

In the following sections, let $\A$ be a Hom-finite abelian category over a field $k$. Assume that $\A$ has enough projectives $\mathcal P$.

\begin{lem}
Any triangle $Q\xrightarrow{q} S\xrightarrow{s} T\xrightarrow{t} Q[1]$ in $ \mathrm{D}^b(\A)$ where $Q,S,T\in\h$ induces an exact sequence $Q\xrightarrow{\overline q} S\xrightarrow{\overline s} T$ in $\overline \h$.
%$\mathbb{F}(R)\xrightarrow{\mathbb{F}(r)} \mathbb{F}(S)\xrightarrow{\mathbb{F}(s)} \mathbb{F}(T)\to 0$ in $\mod \R$.
\end{lem}

\begin{proof}
It is enough to show $\mathbb{F}(R)\xrightarrow{\mathbb{F}(r)} \mathbb{F}(S)\xrightarrow{\mathbb{F}(s)} \mathbb{F}(T)$ is exact in $\mod \mathcal P$. But this is immediately followed by the basic property of triangulated category.
%Let $R'$ be an arbitrary object in $\R$ and $r':R'\to T$ be any morphism. Since $tr'=0$, $r'$ factors through $s$, Hence $\mathbb{F}(s)$ is an epimorpshim. It is obvious that $\mathbb{F}(s)\circ\mathbb{F}(r)=0$, if we have a morphism $r'':R'\to S$ such that $r''s=0$, then $r''$ factors through $r$. Hence we have $\Im \mathbb{F}(r)=\Ker \mathbb{F}(s)$.
\end{proof}

%\begin{prop}
%If $\overline \h$ is abelian, then $\overline \Y=\mathrm{P}(\Fac \overline \U)$ is a support $\tau$-tilting subcategory.
%\end{prop}

We have the following theorem.

\begin{thm}\label{main1}
Let $\A$ be a Hom-finite abelian category with enough projectives $\mathcal P$. Then any covariantly finite $\tau$-rigid subcategory $\V$ is contained in a support $\tau$-tilting subcategory $\mathbf{P}(\Fac \V)$, which is just the following subcategory
$$\add \{ W \in \A \text{ }|\text{ }  \text{there is an exact sequence } P\xrightarrow{v} V\to W\to 0, P\in \mathcal P, v \text{ is a left } \V\text{-approximation} \}.$$
\end{thm}

\begin{proof}
Since $\V$ is covariantly finite, any projective object $P$ admits an exact sequence $P\xrightarrow{v} V\to W\to 0$ where $v$ is a left $\V$-approximation. We show that $W\in \mathbf{P}(\Fac \V)$.

Obviously $W\in \Fac \V$.  Let $P\xrightarrow{ v_1} U\xrightarrow{ v_2} V$ be an epic-monic factorization of $v$. Then we have a short exact sequence $0\to U\xrightarrow{v_2} V\to W\to 0$. %Since $p$ is a left $\U$-approximation, we get $\overline p_2$ is a left $\overline \U$-approximation.
Now let $Z\in \Fac \V$. It admits an epimorphism $z:V'\to Z\to 0$ where $V'\in \V$. Let $u:U\to Z$ be any morphism, since $P$ is projective, there is a morphism $p':P\to V'$ such that $zp'= uv_1$. Since $v$ is a left $\V$-approximation, there is a morphism $v':P\to U'$ such that $v'v=p'$. Hence we have $u=zv'v_2$. We have the following exact sequence
$$\Hom_{\A}(V,Z)\xrightarrow{\Hom_{\A}( v_2,Z)} \Hom_{\A}(U,Z)\xrightarrow{0} \Ext^1_{\A}(W,Z) \to \Ext^1_{\A}(V,Z)=0$$
which implies $\Ext^1_{\A}(X,Z)=0$. Hence $W\in {\rm\bf P}(\Fac \V)$.

By Lemma \ref{imp}, $\mathbf{P}(\Fac \V)$ is a $\tau$-rigid subcategory.  Denote $$\add \{ W \in \A \text{ }|\text{ }  \text{there is an exact sequence } P\xrightarrow{v} V\to W\to 0, P\in \mathcal P, v \text{ is a left } \V\text{-approximation} \}$$ by $\M$, since $\M\subseteq \mathbf{P}(\Fac \V)$, $\M$ is $\tau$-rigid and then by definition it is a support $\tau$-tilting subcategory.
Since $\xymatrix{\A \ar[r]_-{\simeq}^-{\mathbb{F}|_{\A}} &\mod \mathcal P}$, we have $\mathbb{F}(\mathbf{P}(\Fac \V))$ is a $\tau$-rigid subcategory in $\mod \mathcal P$. By Theorem \ref{ZZ}, there is a two-term $\mathcal P[1]$-rigid subcategory $\s$ such that $\mathbb{F}(\s)=\mathbb{F}(\mathbf{P}(\Fac \V))$.  We also have $\mathbb{F}(\M)$ is a support $\tau$-tilting subcategory in $\mod \mathcal P$. By Theorem \ref{ZZ}, there is a two-term weak $\mathcal P[1]$-cluster tilting subcategory $\N$ such that $\mathbb{F}(\M)=\mathbb{F}(\N)$. For a projective object $P$, we have $\Hom_{\A}(P, \M)=0$ if and only if $\Hom_{\A}(P, \mathbf{P}(\Fac \V))=0$, hence $\N$ is contained in $\s$, which implies $\N=\s$ and then we have $\mathbf{P}(\Fac \V)=\M$.
\end{proof}

By the proof of this theorem, we get the following corollary.

\begin{cor}\label{support}
Let $\A$ be a Hom-finite abelian category with enough projectives $\mathcal P$. If $\M$ is a support $\tau$-tilting subcategory, then $\mathrm{P}(\Fac \M)=\M$.
\end{cor}

Moreover, we also have the following theorem.

\begin{thm}\label{main2}
Let $\A$ be an abelian category with enough projectives $\mathcal P$. There is a bijection $\Phi: \M\mapsto\Fac\M$ from the first of the following sets to the second:
\begin{itemize}
\item[(a)] Support $\tau$-tilting subcategories $\M$.
\smallskip

\item[(b)] Subcategories $\s=\Fac \mathbf{P}(\s)$ such that every projective object $P$ admits a left $ \mathbf{P}(\s)$-approximation.
\end{itemize}
\end{thm}

\begin{proof}
We show that $\Phi: \M\mapsto \Fac \M$ is a bijection from the set (a) to the set (b).

By Theorem \ref{main1} and Corollary \ref{support} we know that $\Phi$ is injective.
Now we show that $\Phi$ is surjective.
\smallskip

Let $\mathbf{P}(\s)=\M$. By definition we have $\Ext^1_{\A}(\M,\Fac \M)=0$. By Lemma \ref{imp}, $\M$ is $\tau$-rigid. Let $P$ be any projective object. It admits an exact sequence $P\xrightarrow{p} M^0\to M^1\to 0$ where $p$ is a left $\M$-approximation. We have $M^1\in \Fac \M= \s$. Let $P \xrightarrow{ p_1} Q\xrightarrow{ p_2} M^0$ be an epic-monic factorization of $p$. Then we have a short exact sequence $0\to Q\xrightarrow{p_2} M^0\to M^1\to 0$ where $p_2$ is a left $\M$-approximation. For any object $S\in \s$, we have an exact sequence
$$\Hom_{\A}(M^0,S)\xrightarrow{\Hom_{\A}(p_2,S)} \Hom_{\A}(Q,S)\xrightarrow{0} \Ext^1_{\A}(M^1,S) \to \Ext^1_{\A}(M^0,S)=0$$
which implies that $ \Ext^1_{\A}(M^1,S)=0$. Hence $M^1\in \M$ and by the definition $\M$ is a support $\tau$-tilting subcategory.
\end{proof}

\begin{rem}
The subcategory $\s=\Fac \mathbf{P}(\s)$ is a finitely generated torsion class since $\mathbf{P}(\s)$ is a $\tau$-rigid subcategory.
\end{rem}

%Since $V$ is functorially finite, it is also contained in another two-term weak $\mathcal P[1]$-
%cluster tilting subcategories:
%$$\X=\add \{ X \in \Fac V \text{ }|\text{ }  \text{there is a triangle } P\xrightarrow{v} V\to X\to P[1], v \text{ is a left } (\add V)\text{-approximation} \}.$$

Let $U$ be an object in $\A$ such that $\add U$ is $\tau$-rigid. Then we can get a two-term $\mathcal P[1]$-rigid object $H^{-1}\circ \mathbb{F}(U)=:V$ (we can assume that $V$ does not have any direct summand in $\mathcal P[1]$). Now by \cite[Theorem 3.1]{ZZ}, $V$ is contained in a two-term weak $\mathcal P[1]$-cluster tilting subcategory
$$\Y=\add \{ Y\in \mathrm{D}^b(\A) \text{ }|\text{ } \exists \text{ a triangle } P\xrightarrow{p} Y\to V_0\xrightarrow{v} P[1], P\in \mathcal P,v \text{ is a right } (\add V)\text{-approximation} \}. $$
Since $V$ does not contain any direct summand in $\mathcal P[1]$, we get $p\neq 0$. Hence we have a $\tau$-tilting subcategory $\N=(\mathbb{F}|_{\A})^{-1}\mathbb{F}( \Y)$ that contains $U$. We call $\N$ the \emph{Bongartz completion} of $U$. Moreover, we have the following proposition.

\begin{prop}
If $\A=\mod \Lambda$ where $\Lambda$ is a finite dimensional $k$-algebra, then the Bongartz completion of $U$ is just $\add \mathbf{P}({^{\bot}}(\tau U))$ {\rm (}$\mathbf{P}({^{\bot}}(\tau U))$ is the Bongartz completion of $U$ defined in \cite{AIR}{\rm )}, where ${^{\bot}}(\tau U)=\{X\in \mod \Lambda \text{ }|\text{ }\Hom_{\Lambda}(X,\tau U)=0 \}$.
\end{prop}

\begin{proof}
By \cite[Theorem 2.10]{AIR}, $\mathbf{P}({^{\bot}}(\tau U))$ is a $\tau$-tilting module. Then we have a two-term weak $\Lambda[1]$-cluster tilting subcategory $\Z\subseteq \mathbf{D}^b(\Lambda)$ such that $\mathbb{F}(\Z)=\mathbb{F}(\add\mathbf{P}({^{\bot}}(\tau U)))$. Then $V\in \Z$ and $[\Lambda[1]](V,\Z[1])=0$. Let $Z\in \overline \Z$ and $\overline z:P\to Z\to 0$ be an epimorphism where $P$ is projective. Since $P$ admits a triangle $P\to Y\to V_0\xrightarrow{v} P[1]$, $v$ is a right $(\add V)$-approximation, we have the following commutative diagram
$$\xymatrix{
V_0[-1] \ar[r] &P\ar[r] \ar[d]^z &Y\ar[r] \ar@{.>}[dl]^y &V_0\\
&Z
}
$$
which implies $Z\in \Fac \overline \Y$. Then $\mathbf{P}({^{\bot}}(\tau U))\in \Fac \N$ and $\Fac \mathbf{P}({^{\bot}}(\tau U))\subseteq \Fac \N$. Since $\Ext^1_{\Lambda}(U, \Fac \N)=0$, we have $\Hom_{\Lambda}(\N, \tau U)=0$. Then by \cite[Theorem 2.10]{AIR}, we have $\N\subseteq {^{\bot}}(\tau U)=\Fac \mathbf{P}({^{\bot}}(\tau U))$. Hence we have $\Fac \N=\Fac \mathbf{P}({^{\bot}}(\tau U))$, which implies $\N=\add \mathbf{P}({^{\bot}}(\tau U))$.
\end{proof}

We have the following corollary immediately.

\begin{cor}
Let $\Lambda$ be a finite dimensional $k$-algebra and $U\in \mod\Lambda$ be a $\tau$-rigid module. Then any projective module $P$ admits an exact sequence $P\to Y_0\to U_0\to 0$ where $Y_0\in \add \mathbf{P}({^{\bot}}(\tau U))$ and $U_0\in \add U$.
\end{cor}

%\begin{prop}
%Let $U$ be an object in $\A$ such that $\add U$ is $\tau$-rigid. %Then we can get a two-term $\mathcal P[1]$-rigid object $H^{-1}\circ \mathbb{F}(U)=:V$ (we can assume $V$ does not have any direct summand in $\mathcal P[1]$).  $V$ is contained in two different two-term weak $\mathcal P[1]$-cluster tilting subcategories:
%$$\X=\add \{ X \in \Fac V \text{ }|\text{ }  \text{there is a triangle } P\xrightarrow{v} V\to X\to P[1], v \text{ is a left } (\add V)\text{-approximation} \},$$
%$$\Y=\add \{ Y\in \mathrm{D}^b(\A) \text{ }|\text{ } \text{there is a triangle } P\to Y\to V'\xrightarrow{v'} P[1], v' \text{ is a right } (\add V)\text{-approximation} \}. $$
%Assume $U$ is not support $\tau$-tilting, for any object $X\in \X$ such that $X\notin \add U$,
%\end{prop}

\section{Applications}

In this section, we also assume that $\A$ has enough injectives $\mathcal I$.

\begin{defn}\label{tilting}
Let $\C$ be a subcategory of $\A$.
\begin{itemize}
\item[(i)] $\C$ is said to be hereditary if for any object $C\in \C$, we have $\Ext^2_{\A}(C,-)=0$.
\item[(ii)] $\C$ is said to be partial tilting if $\C$ is hereditary and $\Ext^1_{\A}(\C,\C)=0$.
\item[(iii)] $\C$ is said to be tilting if it is partial tilting and any projective object $P$ admits a short exact sequence
$$0\to P\to C^0\to C^1\to 0$$
where $C^0,C^1\in \C$.
\end{itemize}
\end{defn}

\begin{rem}
By Lemma \ref{imp}, we know that any partial tilting subcategory is $\tau$-rigid.
\end{rem}

%\subsection{Bongartz completion for partial tilting subcategories}

For any subcategory $\C\subseteq \A$, let $\C^{\bot_1}=\{A\in \A \text{ }|\text{ }\Ext^1_{\A}(\C,A)=0 \}$. We have the following lemma.

\begin{prop}\label{completion2}
Let $\U$ be a contravariantly finite partial tilting subcategory. Then $\mathbf{P}(\U^{\bot_1})$ is a tilting subcategory that contains $\U$.
\end{prop}

\begin{proof}
Let $$\M=\add \{M\in  \mathbf{P}(\U^{\bot_1}) \text{ }| \text{ }  \text{there is a short exact sequence } 0\to P\xrightarrow{p} M\to U\to 0, P\in \mathcal P, U\in \mathcal U \}.$$
Since we have
$$0=\Ext^2_{\A}(U,\A)\to \Ext^2_{\A}(M,\A)\to \Ext^2_{\A}(P,\A)=0,$$
$\M$ is a partial tilting subcategory. For any object $P\in \mathcal P$, since we assume $\A$ has enough injectives, $P$ admits a short exact sequence $0\to P\to I\to Q\to 0$ where $I\in \mathcal I$. Let $u\colon U\to Q$ be a right $\U$-approximation. Then we have the following commutative diagram
$$\xymatrix{
0 \ar[r] &P \ar[r] \ar@{=}[d] &M' \ar[r] \ar[d] &U \ar[r] \ar[d]^u &0\\
0 \ar[r] &P \ar[r] &I \ar[r] &Q \ar[r] &0
}
$$
where $M'\in \U^{\bot_1}$.  For any object $N\in \U^{\bot_1}$, we have the following exact sequence
$$0=\Ext^1_{\A}(U,N)\to \Ext^1_{\A}(M',N)\to \Ext^1_{\A}(P,N)=0.$$
Then $M'\in \M$. Hence by the definition, $\M$ is a tilting subcategory of $\A$. Obviously it is a $\tau$-tilting subcategory, so is $\mathbf{P}(\U^{\bot_1})$. By Theorem \ref{main1} and Theorem \ref{ZZ}, the two-term weak $\mathcal P[1]$-cluster tilting subcategory which is correspondent to $\M$ is contained in the one that is
correspondent to $\mathbf{P}(\U^{\bot_1})$, which means they are the same, then by the one-to-one correspondence, we
have $\M=\mathbf{P}(\U^{\bot_1})$.
%Since $\U$ is hereditary, we have $\Fac (\U^{\bot_1})=\U^{\bot_1}$. Then by Lemma \ref{imp} and definition, $\mathbf{P}(\U^{\bot_1})$ is a support $\tau$-tiling subcategory that contains $\M$. By Theorem \ref{main1} and \cite[Theorem 4.5]{ZZ}, the two-term weak $\mathcal P[1]$-cluster tilting subcategory which is correspondent to $\M$ is contained in the one that is correspondent to $\mathbf{P}(\U^{\bot_1})$, which means they are the same, then by the one-to-one correspondence, we have $\M=\mathbf{P}(\U^{\bot_1})$.
\end{proof}

%Let $U$ be an object in $\A$ such that $\add U$ is $\tau$-rigid. Then we can get a two-term $\mathcal P[1]$-rigid object $H^{-1}\circ \mathbb{F}(U)=:V$ (we can assume $V$ does not have any direct summand in $\mathcal P[1]$). Now by \cite[Theorem 3.1]{ZZ}, $V$ is contained in a two-term weak $\mathcal P[1]$-cluster tilting subcategory
%$$\Y=\add \{ Y\in \mathrm{D}^b(\A) \text{ }|\text{ } \exists \text{ a triangle } P\xrightarrow{p} Y\to V_0\xrightarrow{v} P[1], P\in \mathcal P,v \text{ is a right } (\add V)\text{-approximation} \}. $$
%Since $V$ does not contain any direct summand in $\mathcal P[1]$, we get $p\neq 0$. Hence we have a $\tau$-tilting subcategory $\N=(\mathbb{F}|_{\A})^{-1}\mathbb{F}( \Y)$ that contains $U$. We call $\N$ the \emph{Bongartz completion} of $U$. Moreover, we have the following proposition.

%\begin{prop}
%Let $\U=\add U$ be a partial tilting subcategory. Then $\mathbf{P}(\U^{\bot_1})=\N$.
%\end{prop}

\begin{prop}\label{eq}
Let $\U$ be a contravariantly finite partial tilting subcategory. The following conditions are equivalent:
\begin{itemize}
\item[(a)] $\U$ is a tilting subcategory.
\smallskip

\item[(b)] $\U^{\bot_1}=\Fac \U$.
\end{itemize}
\end{prop}

\begin{proof}
By Lemma \ref{imp}, $\U^{\bot_1}\supseteq \Fac \U$ holds. To show (a) implies (b), we only need to prove $\U^{\bot_1}\subseteq \Fac \U$.

Let $X\in \U^{\bot_1}$. It admits an epimorphism $P\xrightarrow{p} X\to 0$ where $P$ is projective. Since $\U$ is tilting, $P$ admits a short exact sequence $0\to P\xrightarrow{m} U^0\to U^1\to 0$ where $U^0,U^1\in \U$. Hence there is a morphism $m'\colon M\to X$ such that $p=m'm$. Moreover, $m'$ is an epimorphism, which implies $X\in \Fac \U$.

Now we show that (b) implies (a).

If $\U$ is not tilting, then by Proposition \ref{completion2} we have a tilting subcategory $\mathbf{P}(\U^{\bot_1})=\mathbf{P}(\Fac \U)\supsetneq \U$. Let $X$ be an indecomposable object in $\mathbf{P}(\U^{\bot_1})$ such that $X\notin \U$. Since $X\in \Fac \U$ and $\U$ is contravariantly finite, $X$ admits a short exact sequence $0\to Y\to U_X\xrightarrow{g} X\to 0$ where $g$ is a right $\U$-approximation. Then $Y\in \U^{\bot_1}=\Fac \U$ and this short exact sequence splits, which implies that $X$ is a direct summand of $U$, a contradiction. Hence $\U$ has to be tilting.
\end{proof}

%\section{mutations of support $\tau$-tilting subcategories}

%When $\A$ is hereditary, any subcategory $\U\subseteq \A$ is $\tau$-rigid if and only if it is partial tilting, if and only if it is rigid. We have the following proposition.

When $\A$ is hereditary, any subcategory $\U\subseteq \A$ is $\tau$-rigid if and only if it is partial tilting, if and only if $\Ext^1_{\A}(\U,\U)=0$. We have the following proposition.

\begin{prop}\label{til}
Let $\A$ be hereditary and $\U\subseteq \A$ be a functorially finite $\tau$-rigid subcategory. Assume that $\U$ is not support $\tau$-tilting, then $\U$ is contained in at least two different tilting subcategories.
\end{prop}

\begin{proof}
By Theorem \ref{main1} and Proposition \ref{completion2}, $\U$ is contained in a tilting subcategories $\mathbf{P}(\U^{\bot_1})=:\Y$ and a support $\tau$-tilting subcategory $\mathbf{P}(\Fac \U)=:\X$.

Let $X\in \X$ be an indecomposable object such that $X\notin \U$. Since $X\in \Fac \U$ and $\U$ is functorially finite, we have a short exact sequence $0\to Y\to U_0\xrightarrow{u} X\to 0$ where $u$ is a right $\U$-approximation, then $Y\in \U^{\bot_1}$. Since $\A$ is hereditary, we have the following exact sequence $$0=\Ext^1_{\A}(U_0, \U^{\bot_1}) \to \Ext^1_{\A}(Y, \U^{\bot_1})\to \Ext^2_{\A}(X, \U^{\bot_1})=0$$ which implies that $Y\in \Y$. Obviously $X\notin \Y$ and $Y\notin \X$, otherwise the short exact sequence splits and $X$ becomes a direct summand of $U_0$, a contradiction.

Let $\U'=\add(\U\cup\{X\})$. Then $\U'$ is a functorially finite $\tau$-rigid subcategory. By Proposition \ref{completion2}, $\U'$ is contained in a tilting subcategories $\mathbf{P}({\U'}^{\bot_1})=:\Y'$. We have $Y\notin \Y'$ since $\Ext^1_{\A}(X,Y)\neq 0$. Hence $\Y$ and $\Y'$ are two different tilting subcategories that contain $\U$.
\end{proof}

\begin{defn}
A covariantly finite $\tau$-rigid subcategory $\U\subseteq \A$ is called almost support $\tau$-tilting if for any support $\tau$-tilting subcategory $\X\supseteq \U$, there exists an indecomposable object $X\notin \U$ such that $\X=\add(\U\cup\{ X\})$.
\end{defn}
%\section{mutations of support $\tau$-tilting subcategories}

%When $\A$ is hereditary, any subcategory $\U\subseteq \A$ is $\tau$-rigid if and only if it is partial tilting, if and only if it is rigid. We have the following theorem.

\begin{thm}\label{main3}
Let $\A$ be hereditary and $\U\subseteq \A$ be a functorially finite almost support $\tau$-tilting subcategory. Then $\U$ is contained in exactly two different tilting subcategories $\mathbf{P}(\U^{\bot_1})$ and $\mathbf{P}(\Fac \U)$.
\end{thm}

To show this theorem, we need the following lemma, which is an analogue of \cite[Lemma 2.20]{AIR}.

\begin{lem}\label{imp2}
Let $\T$ be a contravariantly finite tilting subcategory and $\U$ be a $\tau$-rigid subcategory such that $\T^{\bot_1}\subseteq \U^{\bot_1}$, then any object $U\in \U$ admits a short exact sequence $0\to U\xrightarrow{f} T\xrightarrow{g} T_0\to 0$ where $f$ is a minimal left $(\Fac \T)$-approximation, $T,T_0\in \T$ and $\add T\cap \add T_0=0$.
\end{lem}

\begin{proof}
Let $U\xrightarrow{f} T$ be a minimal left $\T$-approximation approximation. We show that $f$ is a left $(\Fac \T)$-approximation approximation and it is a monomorphism.\\
Let $u:U\to M$ be any morphism where $M\in \Fac \T$. Since $\T$ is a contravariantly finite, $M$ admits a short exact sequence $0\to N\to T'\xrightarrow{t'} M\to 0$ where $t'$ is a right $\T$-approximation. Then $N\in \T^{\bot_1}\subseteq \U^{\bot_1}$. Hence there is a morphism $t:U\to T'$ such that $u=t't$. Then $t$ factors through $f$, which implies that $u$ factors through $f$. Hence $f$ is a left $(\Fac \T)$-approximation approximation. Since any injective object $I\in \T^{\bot_1}=\Fac \T$ by Proposition \ref{eq} and $\A$ has enough injectives, we get that $f$ is a monomorphism.\\
Now we have a short exact sequence $0\to U\xrightarrow{f} T\xrightarrow{g} T_0\to 0$, hence $T_0\in \Fac \T=\T^{\bot_1}$. Since $f$ is a left $\T^{\bot_1}$-approximation and $\Ext^1_{\A}(T,\T^{\bot_1})=0$, we have $T_0\in \mathbf{P}(\T^{\bot_1})=\T$.\\
The proof of $\add T\cap \add T_0=0$ is just the same as in \cite[Lemma 2.20]{AIR}.
\end{proof}

Now we are ready to show Theorem \ref{main3}.

\begin{proof}
We first show that any support $\tau$-tilting subcategory $\V\supsetneq \U$ is tilting.\\
Since $\U$ is almost support $\tau$-tilting, there is an indecomposable $V\notin \U$ such that $\V=\add(\U\cup\{ V\})$. Hence $\V$ is a contravariantly finite $\tau$-rigid subcategory. If $\V$ is not tilting, by Proposition \ref{completion2}, we have a tilting subcategory $\mathbf{P}(\V^{\bot_1})\supsetneq \V$. Since any tilting subcategory is support $\tau$-tilting, we also have an indecomposable object $T\notin \V$ such that $\mathbf{P}(\V^{\bot_1})=\add(\U\cup\{T\})$. But then $V\simeq T$, a contradiction. Hence $\V$ is tilting.

By Theorem \ref{main1} and Proposition \ref{completion2}, $\U$ is contained in two tilting subcategories $\mathbf{P}(\U^{\bot_1})$ and $\mathbf{P}(\Fac \U)$. As in the proof of Proposition \ref{til}, we can get $\mathbf{P}(\U^{\bot_1})\nsupseteq \mathbf{P}(\Fac \U)$ and $\mathbf{P}(\U^{\bot_1})\nsubseteq \mathbf{P}(\Fac \U)$.

%Let $\mathbf{P}(\Fac \U)=\add(\U\cup\{ X\})$ where $X\notin \U$ is indecomposable. Since $X\in \Fac \U$ and $\U$ is functorially finite, we have a short exact sequence $0\to Y_0\to U_0\xrightarrow{u} X\to 0$ where $u$ is a right $\U$-approximation, then $Y_0\in \U^{\bot_1}$. Since $\A$ is hereditary, we have the following exact sequence $0=\Ext^1_{\A}(U_0, \U^{\bot_1}) \to \Ext^1_{\A}(Y_0, \U^{\bot_1})\to \Ext^2_{\A}(X, \U^{\bot_1})=0$ which implies that $Y\in \mathbf{P}(\U^{\bot_1})$. Obviously $X\notin \mathbf{P}(\U^{\bot_1})$ and $Y\notin \mathbf{P}(\Fac \U)$, otherwise the short exact sequence splits and $X$ becomes a direct summand of $U$, a contradiction.\\
Now let $\Z$ be a tilting subcategory that contains $\U$. Then there is an indecomposable $Z\notin \U$ such that $\Z=\add(\U\cup\{ Z\})$. We show either $Z\in \mathbf{P}(\Fac \U)$ or $Z\in \mathbf{P}(\U^{\bot_1})$.

Let $\mathbf{P}(\U^{\bot_1})=\add(\U\cup\{ Y\})$. By Lemma \ref{imp2}, there is a short exact sequence
$$0\to Y\oplus U \xrightarrow{\left({\begin{smallmatrix}
f&\\
&1
\end{smallmatrix}}\right)} Z'\oplus U \to Z''\to 0$$
where $U$ is an arbitrary object in $\U$, $\left({\begin{smallmatrix}
f&\\
&1
\end{smallmatrix}}\right)$ is a minimal left $(\Fac \Z)$-approximation. $Z',Z''\in \Z$ and $\add (Z'\oplus U)\cap \add Z''=0$. Then $Z''\in \add Z$.

If $Z''\neq 0$, then we have $Z'\in \U$ and $Z''\in \Fac \U$. Hence $Z\in \Fac \U$. Since $\Ext^1_{\A}(Z,\Fac \U)=0$, we have $Z\in \mathbf{P}(\Fac \U)$.

If $Z''=0$, then $Z'\simeq Y$. Hence $Z\simeq Y$ and $Z\in \mathbf{P}(\U^{\bot_1})$.
\end{proof}

%\begin{thm}
%Let $\U$ be a functorially finite $\tau$-rigid subcategory of $\A$. Assume $\U$ is not support $\tau$-tilting. Then $\U$ is contained in a support $\tau$-titling subcategory $\Z$ such that $\Z\nsubseteq \mathbf{P}(\Fac \U)$ and $\Z\nsupseteq \mathbf{P}(\Fac \U)$.
%\end{thm}

%\begin{proof}

%\end{proof}

\end{document}